\documentclass{article}
\usepackage{amsfonts}
\usepackage{amsmath}
\usepackage[doublespacing]{setspace}

\newtheorem{theorem}{Theorem}
\newtheorem{acknowledgement}[theorem]{Acknowledgement}

\newtheorem{example}[theorem]{Example}

\input{tcilatex}
\begin{document}

\title{SURFACE PENCIL WITH A COMMON LINE OF CURVATURE IN MINKOWSKI 3-SPACE}
\author{Evren ERG\"{U}N, Ergin BAYRAM, \ Emin KASAP \\
Ondokuz May\i s University, Faculty of Arts and Sciences, \\
Department of Mathematics, Samsun TURKEY\\
eergun@omu.edu.tr, \ erginbayram@yahoo.com, \ kasape@omu.edu.tr}
\date{}
\maketitle

\begin{abstract}
In this paper, we analyze the problem of constructing a surface pencil from
a given spacelike (timelike) line of curvature. By using the Frenet frame of
the given curve in Minkowski 3-space, we express the surface pencil as a
linear combination of this frame and derive the necessary and sufficient
conditions for the coefficients to satisfy the line of curvature
requirement. To illustrate the method some examples showing members of the
surface pencil with their line of curvature are given.
\end{abstract}

\textbf{Keywords:} Minkowski space, Line of curvature, Surface pencil.

\textbf{2000} \textbf{Mathematics Subject Classification:} 53B30, 51B20,
53C50.

\textbf{Short Title:} Surfaces with common line of curvature in $%
\mathbb{R}
_{1}^{3}$

\bigskip

\bigskip

\bigskip

\section{Introduction}

On a Minkowski surface, tangent vectors are classified into timelike,
spacelike or null and so a curve on the surface is said to be timelike,
spacelike or null if its tangent vectors are always timelike, spacelike or
null, respectively. In fact, a timelike curve corresponds to the path of an
observer moving slower than the speed of light, spacelike curve corresponds
to moving faster than the speed of light and null curve corresponds to
moving at the speed of light. Surface with a common characteristic curve
have been the subject of many recent studies. Wang et.al.[1] studied the
problem of constructing a surface pencil from a given spatial geodesic. They
parametrized the surface by using the Frenet frame of the given curve and
gave the necessary and sufficient condition to satify the geodesic
requirement. Kasap and Aky\i ld\i z [2] considered the surfaces with common
geodesic in Minkowski 3-space. They studied spacelike surface with common
spacelike geodesic and timelike surface with common spacelike or timelike
geodesic. Li et.al.[3] derived the necessary and sufficient condition for a
given curve to be a line of curvature on a surface. \c{S}affak and Kasap [4]
analyzed the problem of finding a surface family through a null geodesic
with Cartan frame. Being inspired by the above studies, we extend the method
of Lie et. al. [3] to derive the necessary and sufficient condition for the
given curve to be a line of curvature for the parametric surface. By
utilizing the Frenet frame, we derive necessary and sufficient condition for
the correct parametric representation of the surface $P\left( s,t\right) $
when the parameter $s$ is the arc-length of the curve $r\left( s\right) $
and find the necessary constraints on the coefficients of vectors of the
frame ( which are called marching-scale functions) so that both the line of
curvature and parametric requirement are met. Thus, we defined the spacelike
and timelike surface pencil with common line of curvature. Also, we give two
general forms of the marching-scale functions and obtain necessary and
sufficient conditions on them for which the given curve is a line of
curvature on a surface. Finally, we give some nice examples, showing members
of the surface pencil with their common line of curvatures to illustrate the
method.

\section{Preliminaries}

Let us consider Minkowski 3-space $%
\mathbb{R}
_{1}^{3}=\left[ 
\mathbb{R}
^{3},\left( +,+,-\right) \right] $ and let the Lorentzian inner product of $%
X=\left( x_{1},x_{2},x_{3}\right) $ and $Y=\left( y_{1},y_{2},y_{3}\right)
\in 
\mathbb{R}
_{1}^{3}$ be

\begin{equation*}
\left\langle X,Y\right\rangle =x_{1}y_{1}+x_{2}y_{2}-x_{3}y_{3}.
\end{equation*}%
A vector $X\in 
\mathbb{R}
_{1}^{3}$ is called a spacelike vector when $\left\langle X,X\right\rangle
>0 $ or $X=0.$ It is called timelike and null (lightlike) vector in case of $%
\left\langle X,X\right\rangle <0$, $\left\langle X,X\right\rangle =0$ for $%
X\neq 0.$respectively, [5]. The vector product of vectors $X=\left(
x_{1},x_{2},x_{3}\right) $ and $Y=\left( y_{1},y_{2},y_{3}\right) \in 
\mathbb{R}
_{1}^{3}$ is defined by [6]

\begin{equation*}
X\times Y=\left(
x_{2}y_{3}-x_{3}y_{2},~x_{3}y_{1}-x_{1}y_{3},x_{2}y_{1}-x_{1}y_{2}\right) .
\end{equation*}

Let $\ r=r\left( s\right) $ be a unit speed curve in $%
\mathbb{R}
_{1}^{3}$. By $\kappa \left( s\right) $ and $\tau \left( s\right) $ we
denote the natural curvature and torsion of $r\left( s\right) $
,respectively. Consider the Frenet frame$\left\{ T,N,B\right\} $ associated
with the curve $\ r=r\left( s\right) $ such that $T=T\left( s\right)
,N=N\left( s\right) $ and $B=B\left( s\right) $ are the unit tangent, the
principal normal and the binormal vector fields, respectively. If $r=r\left(
s\right) $ is a spacelike curve, then the structural equations (or Frenet
formulas) of this frame are given as%
\begin{equation*}
\overset{\cdot }{T}\left( s\right) =\kappa \left( s\right) N\left( s\right)
,~\overset{\cdot }{N}\left( s\right) =\varepsilon \kappa \left( s\right)
T\left( s\right) +\tau \left( s\right) B\left( s\right) ,~\overset{\cdot }{B}%
\left( s\right) =\tau \left( s\right) N\left( s\right)
\end{equation*}%
where 
\begin{equation*}
\varepsilon =\{%
\begin{array}{c}
+1,B~is~spacelike, \\ 
-1,~B~is~timelike.%
\end{array}%
\end{equation*}

If $r=r\left( s\right) $ is a timelike curve, then above equations are given
as [7]%
\begin{equation*}
\overset{\cdot }{T}\left( s\right) =\kappa \left( s\right) N\left( s\right)
,~\overset{\cdot }{N}\left( s\right) =\kappa \left( s\right) T\left(
s\right) -\tau \left( s\right) B\left( s\right) ,~\overset{\cdot }{B}\left(
s\right) =\tau \left( s\right) N\left( s\right) .
\end{equation*}%
The norm of a vector $X$ is defined by [5]%
\begin{equation*}
\left\Vert X\right\Vert _{IL}=\sqrt{\left\vert \left\langle X,X\right\rangle
\right\vert }\text{.}
\end{equation*}

\begin{theorem}
: Let $\ X$ and $Y$ be non-zero orthogonal Lorentz vectors in $%
\mathbb{R}
_{1}^{3}.$ If $X$ is timelike, then $Y$ is spacelike \textit{[8]}.
\end{theorem}

\begin{theorem}
: Let $X$ and $Y$ be positive (negative ) timelike vectors in $%
\mathbb{R}
_{1}^{3}$. Then 
\begin{equation*}
\left\langle X,Y\right\rangle \leq \left\Vert X\right\Vert \left\Vert
Y\right\Vert
\end{equation*}%
with equality if and only if $X$ and $Y$ are linearly dependent [8].
\end{theorem}

Let $X$ and $Y$ be positive (negative ) timelike vectors in $%
\mathbb{R}
_{1}^{3}$. Then there is a unique non-negative real number $\varphi \left(
X,Y\right) $ such that%
\begin{equation*}
\left\langle X,Y\right\rangle =\left\Vert X\right\Vert \left\Vert
Y\right\Vert \cosh \varphi \left( X,Y\right)
\end{equation*}%
the Lorentzian timelike angle between $X$ and $Y$ is defined to be $\varphi
\left( X,Y\right) $ [8].

Let $X$\ and $Y$\ be spacelike vectors in $%
\mathbb{R}
_{1}^{3}$\ that span a spacelike vector subspace. Then we have 
\begin{equation*}
\left\vert \left\langle X,Y\right\rangle \right\vert \leq \left\Vert
X\right\Vert \left\Vert Y\right\Vert
\end{equation*}%
with equality if and only if $X$ and $Y$ are linearly dependent. Hence,
there is a unique real number $\varphi \left( X,Y\right) $\ between $0$\ and 
$\pi $\ such that 
\begin{equation*}
\left\langle X,Y\right\rangle =\left\Vert X\right\Vert \left\Vert
Y\right\Vert \cos \varphi \left( X,Y\right)
\end{equation*}%
the Lorentzian spacelike angle between $X$\ and $Y$\ is defined to be $%
\varphi \left( X,Y\right) $ [8].

Let $X$\ and $Y$\ be spacelike vectors in $%
\mathbb{R}
_{1}^{3}$\ that span a timelike vector subspace. Then we have 
\begin{equation*}
\left\vert \left\langle X,Y\right\rangle \right\vert >\left\Vert
X\right\Vert \left\Vert Y\right\Vert .
\end{equation*}%
Hence, there is a unique positive real number $\varphi \left( X,Y\right) $\
between $0$\ and $\pi $\ such that 
\begin{equation*}
\left\vert \left\langle X,Y\right\rangle \right\vert =\left\Vert
X\right\Vert \left\Vert Y\right\Vert \cosh \varphi \left( X,Y\right)
\end{equation*}%
the Lorentzian timelike angle between $X$\ and $Y$\ is defined to be $%
\varphi \left( X,Y\right) $ [8].

Let $X$\ be a spacelike vector and $Y$\ be a positive timelike vector in $%
\mathbb{R}
_{1}^{3}$. Then there is a unique nonnegative real number $\varphi \left(
X,Y\right) $\ such that 
\begin{equation*}
\left\vert \left\langle X,Y\right\rangle \right\vert =\left\Vert
X\right\Vert \left\Vert Y\right\Vert \sinh \varphi \left( X,Y\right)
\end{equation*}%
the Lorentzian timelike angle between $X$\ and $Y$\ is defined to be $%
\varphi \left( X,Y\right) $ [8].

Let $M$ be a semi-Riemannian submanifold of $\overset{\_}{M}$ and $\overset{%
\_}{D}$ be the Levi-Civita connection of $\overset{\_}{M},$the function $%
II:\chi \left( M\right) \times \chi \left( M\right) \longrightarrow \chi
\left( M\right) ^{\perp }$ such that 
\begin{equation*}
II\left( X,Y\right) =nor\overset{\_}{D}_{X}Y
\end{equation*}%
is $\Im \left( M\right) -$bilinear and symmetric. $II$ is callled the shape
tensor (or second fundamental form tensor) of $M\subset \overset{\_}{M}$ [5].

Let $N$ be a unit normal vector field on a semi-Riemannian hypersurface $%
M\subset \overset{\_}{M}$. The $\left( 1,1\right) $ tensor field $S$ on $M$
such that 
\begin{equation*}
\left\langle S\left( X\right) ,Y\right\rangle =\left\langle II\left(
X,Y\right) ,N\right\rangle
\end{equation*}%
for all $X,Y\in \chi \left( M\right) $ is called the shape operator of $%
M\subset \overset{\_}{M}$ derived from $N$ [5].

As usual, $S$ determines a linear operator $S:T_{P}\left( M\right)
\longrightarrow T_{P}\left( M\right) $ at each point $P\in M$. If $S$ is
shape operator derived from $N$, then $S\left( X\right) =-\overset{\_}{D}%
_{X}N$ and at each point the linear operator $S$ on $T_{P}\left( M\right) $
is self adjoint. A regular curve r on $M$ is said to be a line of curvature
of $M$ if for all $p\in r$ the tangent line of $r$ is a principal direction
at $p$. According to this definition, the differential equation of the line
of curvature on $M$ is $S\left( T\right) =\omega T,~\omega \neq 0,$ where $S$
is the shape operator of $M$.

A surface in $%
\mathbb{R}
_{1}^{3}$ is called a timelike surface if the induced metric on the surface
is a Lorentzian metric and is called a spacelike surface if induced metric
on the surface is positive definite Riemannian metric, i.e. the normal
vector on the spacelike (timelike) surface is a timelike (spacelike) vector
[9].

A parametric curve $r\left( s\right) $ is a curve on a surface $P=P\left(
s,t\right) $ in $%
\mathbb{R}
_{1}^{3}$ that has a constant $s$ or $t$ parameter value, that is, there
exists a parameter $s_{0\text{ }}$ or $t_{0\text{ }}$such that 
\begin{equation*}
r\left( s\right) =P\left( s,t_{0}\right) \text{ or }r\left( t\right)
=P\left( s_{0},t\right) .
\end{equation*}

\section{Spacelike Surface Pencil with a Common Line of Curvature}

Let $P=P\left( s,t\right) $ be a parametric spacelike surface and $r=r\left(
s\right) $ be a spacelike curve with spacelike binormal. The surface is
defined by the given curve as%
\begin{equation}
P\left( s,t\right) =r\left( s\right) +\left( u\left( s,t\right) ,v\left(
s,t\right) ,w\left( s,t\right) \right) \left( 
\begin{array}{c}
T\left( s\right) \\ 
N\left( s\right) \\ 
B\left( s\right)%
\end{array}%
\right) ,  \tag{3.1}
\end{equation}%
$L_{1}\leq s\leq L_{2},T_{1}\leq t\leq T_{2},$ where $u\left( s,t\right)
,~v\left( s,t\right) ~$and $w\left( s,t\right) $ are called the
marching-scale functions and $\left\{ T\left( s\right) ,N\left( s\right)
,B\left( s\right) \right\} $ is the Frenet frame associated with the curve $%
r\left( s\right) .$

Since the curve $r\left( s\right) $ is a parametric curve on the surface $%
P\left( s,t\right) $, there exists a parameter $t_{0}\in \left[ T_{1},T_{2}%
\right] $ such that $P\left( s,t_{0}\right) =r\left( s\right) $, $L_{1}\leq
s\leq L_{2},$ that is, 
\begin{equation}
u\left( s,t_{0}\right) =v\left( s,t_{0}\right) =w\left( s,t_{0}\right)
\equiv 0~,~L_{1}\leq s\leq L_{2}.  \tag{3.2}
\end{equation}

Let $n_{1}\left( n_{1}=\cosh \theta N+\sinh \theta B\right) $ be a vector
orthogonal to the curve $r\left( s\right) $, where $\theta =\theta \left(
s\right) $ is the Lorentzian timelike angle between $N$ and $n_{1}.$The
curve $r\left( s\right) $ is a line of curvature on the surface $P\left(
s,t\right) $ if and only if $n_{1}$ is parallel to the normal vector $%
n\left( s,t\right) $ of the surface $P\left( s,t\right) $ and $S\left(
T\right) =\omega T,~\omega \neq 0$, where $S$ is the shape operator of the
surface.

Firstly, we derive the condition for $n_{1}$ to be parallel to the normal
vector $n\left( s,t\right) $ of the surface $P\left( s,t\right) :$

The normal vector can be expressed as

$n\left( s,t\right) =\frac{\partial P\left( s,t\right) }{\partial s}\times 
\frac{\partial P\left( s,t\right) }{\partial t}$

$=(-(\frac{\partial w\left( s,t\right) }{\partial s}+v\left( s,t\right) \tau
\left( s\right) )\frac{\partial v\left( s,t\right) }{\partial t}+$

$(\frac{\partial v\left( s,t\right) }{\partial s}+u\left( s,t\right) \kappa
\left( s\right) +w\left( s,t\right) \tau \left( s\right) )\frac{\partial
w\left( s,t\right) }{\partial t})T\left( s\right) +$

$((1+\frac{\partial u\left( s,t\right) }{\partial s}+v\left( s,t\right)
\kappa \left( s\right) )\frac{\partial w\left( s,t\right) }{\partial t}-$

$(\frac{\partial w\left( s,t\right) }{\partial s}+v\left( s,t\right) \tau
\left( s\right) )\frac{\partial u\left( s,t\right) }{\partial t})N\left(
s\right) +$

$((1+\frac{\partial u\left( s,t\right) }{\partial s}+v\left( s,t\right)
\kappa \left( s\right) )\frac{\partial v\left( s,t\right) }{\partial t}-$

$(\frac{\partial v\left( s,t\right) }{\partial s}+u\left( s,t\right) \kappa
\left( s\right) +w\left( s,t\right) \tau \left( s\right) )\frac{\partial
u\left( s,t\right) }{\partial t})B\left( s\right) $

Thus, we get 
\begin{equation*}
n\left( s,t_{0}\right) =\phi _{1}\left( s,t_{0}\right) T\left( s\right)
+\phi _{2}\left( s,t_{0}\right) N\left( s\right) +\phi _{3}\left(
s,t_{0}\right) B\left( s\right) ,
\end{equation*}%
where

$\phi _{1}\left( s,t_{0}\right) =\frac{\partial v\left( s,t_{0}\right) }{%
\partial s}\frac{\partial w\left( s,t_{0}\right) }{\partial t}-\frac{%
\partial w\left( s,t_{0}\right) }{\partial s}\frac{\partial v\left(
s,t_{0}\right) }{\partial t},$

$\phi _{2}\left( s,t_{0}\right) =\left( 1+\frac{\partial u\left(
s,t_{0}\right) }{\partial s}\right) \frac{\partial w\left( s,t_{0}\right) }{%
\partial t}-\frac{\partial w\left( s,t_{0}\right) }{\partial s}\frac{%
\partial u\left( s,t_{0}\right) }{\partial t},$

$\phi _{3}\left( s,t_{0}\right) =\left( 1+\frac{\partial u\left(
s,t_{0}\right) }{\partial s}\right) \frac{\partial v\left( s,t_{0}\right) }{%
\partial t}-\frac{\partial v\left( s,t_{0}\right) }{\partial s}\frac{%
\partial u\left( s,t_{0}\right) }{\partial t}.$

This follows that $n_{1}\left( s\right) //n\left( s,t_{0}\right) $, $%
L_{1}\leq s\leq L_{2}$, if and only if there exits a function $\lambda
\left( s\right) \neq 0$ such that%
\begin{equation}
\phi _{1}\left( s,t_{0}\right) =0,\phi _{2}\left( s,t_{0}\right) =\lambda
\left( s\right) \cosh \theta ,\phi _{3}\left( s,t_{0}\right) =\lambda \left(
s\right) \sinh \theta .  \tag{3.3}
\end{equation}%
Secondly, since $S\left( T\right) =\omega T,~\omega \neq 0,$ we obtain%
\begin{equation}
\theta \left( s\right) =-\int\limits_{s_{0}}^{s}\tau ds+\theta _{0\text{,}} 
\tag{3.4}
\end{equation}%
where $s_{0}$ is the starting value of the arc-length and $\theta =\theta
\left( s\right) .$ In this paper, we assume that $s_{0}=0.$

Combining (3.2), (3.3) and (3.4), we have the following theorem.

\begin{theorem}
A spacelike curve $r\left( s\right) $ with spacelike binormal is a line of
curvature on the surface $P\left( s,t\right) $ if and only if the followings
are satisfied:
\end{theorem}

$\theta \left( s\right) =-\int\limits_{s_{0}}^{s}\tau ds+\theta \left(
0\right) ,$

$\bigskip u\left( s,t_{0}\right) =v\left( s,t_{0}\right) =w\left(
s,t_{0}\right) \equiv 0,$

$\phi _{1}\left( s,t_{0}\right) \equiv 0~,~\phi _{2}\left( s,t_{0}\right)
=\lambda \left( s\right) \cosh \theta ,~\phi _{3}\left( s,t_{0}\right)
=\lambda \left( s\right) \sinh \theta .$

\bigskip

We call the set of surfaces defined by (3.1) - (3.4) \textit{spacelike
surface pencil with a common line of curvature}. Any surface $P\left(
s,t\right) $ defined by (3.1) and satisfying (3.2) - (3.4) is a member of
this family.

Now, we analyse two different types of the marching-scale functions

$u\left( s,t\right) ,v\left( s,t\right) ~$and$~w\left( s,t\right) $ in. the
Eq. (3.1).

\textbf{(i)} If we choose

$u\left( s,t\right) =\overset{p}{\underset{k=1}{\sum }}a_{1k}l\left(
s\right) ^{k}U\left( t\right) ^{k},$ $v\left( s,t\right) =\overset{p}{%
\underset{k=1}{\sum }}a_{2k}m\left( s\right) ^{k}V\left( t\right) ^{k}$ and$%
\ w\left( s,t\right) =\overset{p}{\underset{k=1}{\sum }}a_{3k}n\left(
s\right) ^{k}W\left( t\right) ^{k}$

then, we can simply express the sufficient condition for which the curve $%
r\left( s\right) $ is a line of curvature of the surface $P\left( s,t\right) 
$ as

\begin{eqnarray}
U\left( t_{0}\right) &=&V\left( t_{0}\right) =W\left( t_{0}\right) =0  \notag
\\
\theta \left( s\right) &=&-\int\limits_{s_{0}}^{s}\tau ds+\theta _{0\text{,}}
\TCItag{3.5} \\
a_{21}m\left( s\right) V^{\prime }\left( t_{0}\right) &=&\lambda \left(
s\right) \sinh \theta ,~a_{31}n\left( s\right) W^{\prime }\left(
t_{0}\right) =\lambda \left( s\right) \cosh \theta ,  \notag
\end{eqnarray}

$\lambda \left( s\right) \neq 0,$ where $l\left( s\right) ,m\left( s\right)
,n\left( s\right) ,U\left( t\right) ,V\left( t\right) $ and $W\left(
t\right) $ are $C^{1}$ functions, $a_{ij}\in 
\mathbb{R}
$ $\left( k=1,2,3;j=1,2,3,...,p\right) $.

\textbf{(ii)} If we choose

$u\left( s,t\right) =f\left( \overset{p}{\underset{k=1}{\sum }}a_{1k}l\left(
s\right) ^{k}U\left( t\right) ^{k}\right) ,$ \ \ $v\left( s,t\right)
=g\left( \overset{p}{\underset{k=1}{\sum }}a_{2k}m\left( s\right)
^{k}V\left( t\right) ^{k}\right) $ \ and $\ w\left( s,t\right) =h\left( 
\overset{p}{\underset{k=1}{\sum }}a_{3k}n\left( s\right) ^{k}W\left(
t\right) ^{k}\right) $

then, we can express the sufficient condition for which the curve $r\left(
s\right) $ is a line of curvature on the surface $P\left( s,t\right) $ as

\begin{eqnarray}
U\left( t_{0}\right) &=&V\left( t_{0}\right) =W\left( t_{0}\right) =0\ and\
f\left( 0\right) =g\left( 0\right) =h\left( 0\right) ,  \notag \\
\theta \left( s\right) &=&-\int\limits_{s_{0}}^{s}\tau ds+\theta _{0\text{,}}
\TCItag{3.6} \\
g^{^{\prime }}\left( 0\right) a_{21}m\left( s\right) V^{\prime }\left(
t_{0}\right) &=&\lambda \left( s\right) \sinh \theta ,~h^{^{\prime }}\left(
0\right) a_{31}n\left( s\right) W^{\prime }\left( t_{0}\right) =\lambda
\left( s\right) \cosh \theta ,  \notag
\end{eqnarray}

$\lambda \left( s\right) \neq 0,$ where $l\left( s\right) ,m\left( s\right)
,n\left( s\right) ,U\left( t\right) ,V\left( t\right) $ and $W\left(
t\right) $ are $C^{1}$ functions, $a_{ij}\in 
\mathbb{R}
$ $\left( k=1,2,3;j=1,2,3,...,p\right) $.

\begin{example}
Let $r\left( s\right) =\left( a\sinh \left( \frac{s}{c}\right) ,\frac{bs}{c}%
,a\cosh \left( \frac{s}{c}\right) \right) $ be a spacelike curve ,
\end{example}

$a,b,c\in 
\mathbb{R}
,~a^{2}+b^{2}=c^{2}~$and $-2\leq s\leq 2.$ It is easy to show that

$T\left( s\right) =\left( \frac{a}{c}\cosh \left( \frac{s}{c}\right) ,\frac{b%
}{c},\frac{a}{c}\sinh \left( \frac{s}{c}\right) \right) ,$

$N\left( s\right) =\left( \sinh \left( \frac{s}{c}\right) ,0,\cosh \left( 
\frac{s}{c}\right) \right) ,$

$B\left( s\right) =\left( \frac{b}{c}\cosh \left( \frac{s}{c}\right) ,-\frac{%
a}{c},\frac{b}{c}\sinh \left( \frac{s}{c}\right) \right) .$

By taking $\theta \left( 0\right) =0$ we have $\theta \left( s\right) =-%
\frac{bs}{c^{2}}$. If we choose $\lambda \left( s\right) \equiv
1,~t_{0}=0,~a_{21}=a_{31}=1$ and

$u\left( s,t\right) =\overset{3}{\underset{k=1}{\sum }}a_{1k}l\left(
s\right) U\left( t\right) \equiv 0,$

$v\left( s,t\right) =\sinh \left( -\frac{bs}{c^{2}}\right) t+\overset{3}{%
\underset{k=2}{\sum }}a_{2k}\sinh ^{k}\left( -\frac{bs}{c^{2}}\right) t^{k},$

$w\left( s,t\right) =\cosh \left( -\frac{bs}{c^{2}}\right) t+\overset{3}{%
\underset{k=2}{\sum }}a_{3k}\cosh ^{k}\left( -\frac{bs}{c^{2}}\right) t^{k}$

then the Eq. (3.5) is satisfied.

Letting $a=b=1$, we immediately obtain a member of the surface pencil (Fig.
3.1) as

$P_{1}\left( s,t\right) =(~\sinh \left( \frac{s}{\sqrt{2}}\right) +\sinh
\left( -\frac{s}{2}\right) t\sinh \left( \frac{s}{\sqrt{2}}\right) +\sinh
\left( \frac{s}{\sqrt{2}}\right) \overset{3}{\underset{k=2}{\sum }}%
a_{2k}\sinh ^{k}\left( -\frac{s}{2}\right) t^{k}+\frac{\sqrt{2}}{2}\cosh
\left( -\frac{s}{2}\right) t\cosh \left( \frac{s}{\sqrt{2}}\right) +\frac{%
\sqrt{2}}{2}\cosh (\frac{s}{\sqrt{2}})\overset{3}{\underset{k=2}{\sum }}%
a_{3k}\cosh ^{k}\left( -\frac{s}{2}\right) t^{k},$

$\frac{s}{\sqrt{2}}-\frac{\sqrt{2}}{2}\cosh \left( -\frac{s}{2}\right) t-%
\frac{\sqrt{2}}{2}\overset{3}{\underset{k=2}{\sum }}a_{3k}\cosh ^{k}\left( -%
\frac{s}{2}\right) t^{k},\cosh (\frac{s}{\sqrt{2}})+\sinh (-\frac{s}{2}%
)t\cosh (\frac{s}{\sqrt{2}})+\cosh (\frac{s}{\sqrt{2}})\overset{3}{\underset{%
k=2}{\sum }}a_{2k}\sinh ^{k}\left( -\frac{s}{2}\right) t^{k}$

$+\frac{\sqrt{2}}{2}\cosh \left( -\frac{s}{2}\right) t\sinh \left( \frac{s}{%
\sqrt{2}}\right) +\frac{\sqrt{2}}{2}\sinh \left( \frac{s}{\sqrt{2}}\right) 
\overset{3}{\underset{k=2}{\sum }}a_{3k}\cosh ^{k}\left( -\frac{s}{2}\right)
t^{k}~).$

\FRAME{dtbpFU}{2.7882in}{2.5633in}{0pt}{\Qcb{Fig. 3.1. $P_{1}\left(
s,t\right) \ $as a member of the surface pencil and its line of curvature}}{%
}{Figure}{\special{language "Scientific Word";type "GRAPHIC";display
"USEDEF";valid_file "T";width 2.7882in;height 2.5633in;depth
0pt;original-width 1.9796in;original-height 1.9796in;cropleft "0";croptop
"1";cropright "1";cropbottom "0";tempfilename
'MDSHIT07.wmf';tempfile-properties "XPR";}}

\begin{example}
Let $r\left( s\right) =\left( \frac{\sqrt{3}}{2}\sinh \left( s\right) ,\frac{%
s}{2},\frac{\sqrt{3}}{2}\cosh \left( s\right) \right) $ be an arc-length
spacelike curve, 0$\leq s\leq 2\pi .$ It is easy to show that
\end{example}

$T\left( s\right) =\left( \frac{\sqrt{3}}{2}\cosh \left( s\right) ,\frac{1}{2%
},\frac{\sqrt{3}}{2}\sinh \left( s\right) \right) ,$

$N\left( s\right) =\left( \sinh \left( s\right) ,0,\cosh \left( s\right)
\right) ,$

$B\left( s\right) =\left( \frac{1}{2}\cosh \left( s\right) ,-\frac{\sqrt{3}}{%
2},\frac{1}{2}\sinh \left( s\right) \right) ,$

$\tau =\frac{1}{2}.$

Taking $\theta \left( 0\right) =0$ we have $\theta \left( s\right) =-\frac{1%
}{2}s$. If we choose $\lambda \left( s\right) \equiv
1,~t_{0}=0,~a_{21}=a_{31}=1$ and

$u\left( s,t\right) =t,~v\left( s,t\right) =\sinh \left( -\frac{1}{2}%
s\right) t,~w\left( s,t\right) =\cosh \left( -\frac{1}{2}s\right) t$

then the Eq. (3.5) is satisfied. So, we have the following surface as a
member of the surface pencil with common line of curvature $r\left( s\right) 
$ (Fig. 3.2) as

$P_{2}\left( s,t\right) =(~\frac{\sqrt{3}}{2}\sinh \left( s\right) +\frac{%
\sqrt{3}}{2}t\cosh \left( s\right) +t\cosh \left( -\frac{1}{2}s\right) \sinh
\left( s\right) -t\sinh \left( -\frac{1}{2}s\right) \cosh \left( -\frac{1}{2}%
s\right) ,$

$\frac{s}{2}+\frac{t}{2}+\frac{\sqrt{3}}{2},$

$\frac{\sqrt{3}}{2}\cosh \left( s\right) +\frac{\sqrt{3}}{2}t\sinh \left(
s\right) +\cosh \left( -\frac{1}{2}s\right) \cosh \left( s\right) -\frac{1}{2%
}\sinh \left( -\frac{1}{2}s\right) \sinh \left( s\right) ~)$

where $0$ $\leq s\leq 2\pi ,~-2\leq t\leq 2$.

\FRAME{dtbpFU}{2.1715in}{2.2295in}{0pt}{\Qcb{Fig. 3.2. $P_{2}\left(
s,t\right) $ as a member of the surface pencil and its line of curvature}}{}{%
Figure}{\special{language "Scientific Word";type
"GRAPHIC";maintain-aspect-ratio TRUE;display "USEDEF";valid_file "T";width
2.1715in;height 2.2295in;depth 0pt;original-width 1.542in;original-height
1.5835in;cropleft "0";croptop "1";cropright "1";cropbottom "0";tempfilename
'MDS86S01.wmf';tempfile-properties "XPR";}}

\bigskip

\section{Timelike Surface Pencil with a Common Line of Curvature}

Let $P=P\left( s,t\right) $ be a parametric timelike surface and $r=r\left(
s\right) $ be a spacelike curve with timelike binormal. The surface is
defined by the given curve as%
\begin{equation}
P\left( s,t\right) =r\left( s\right) +\left( u\left( s,t\right) ,v\left(
s,t\right) ,w\left( s,t\right) \right) \left( 
\begin{array}{c}
T\left( s\right) \\ 
N\left( s\right) \\ 
B\left( s\right)%
\end{array}%
\right) ,  \tag{4.1}
\end{equation}%
$L_{1}\leq s\leq L_{2},T_{1}\leq t\leq T_{2},$ where $\left\{ T\left(
s\right) ,N\left( s\right) ,B\left( s\right) \right\} $ is the Frenet frame
associated with the curve $r\left( s\right) .$

Since the curve $r\left( s\right) $ is a parametric curve on the surface $%
P\left( s,t\right) $, there exists a parameter $t_{0}\in \left[ T_{1},T_{2}%
\right] $ such that $P\left( s,t_{0}\right) =r\left( s\right) $ $L_{1}\leq
s\leq L_{2},$ that is , 
\begin{equation}
u\left( s,t_{0}\right) =v\left( s,t_{0}\right) =w\left( s,t_{0}\right)
\equiv 0~,~L_{1}\leq s\leq L_{2}.  \tag{4.2}
\end{equation}

Let $n_{1}\left( n_{1}=\cosh \theta N+\sinh \theta B\right) $ be a vector
orthogonal to the curve $r\left( s\right) $, where $\theta =\theta \left(
s\right) $ is the Lorentzian timelike angle between $N$ and $n_{1}.$The
curve $r\left( s\right) $ is a line of curvature on the surface $P\left(
s,t\right) $ if and only if $n_{1}$ is parallel to the normal vector $%
n\left( s,t\right) $ of the surface $P\left( s,t\right) $ and $S\left(
T\right) =\omega T,~\omega \neq 0$, where $S$ is the shape operator of the
surface.

Firstly, we derive the condition for $n_{1}$ to be parallel to the normal
vector $n\left( s,t\right) $ of the surface $P\left( s,t\right) :$

The normal vector can be expressed as

$n\left( s,t\right) =\frac{\partial P\left( s,t\right) }{\partial s}\times 
\frac{\partial P\left( s,t\right) }{\partial t}$

$=(-(\tau \left( s\right) v\left( s,t\right) +\frac{\partial w\left(
s,t\right) }{\partial s})\frac{\partial v\left( s,t\right) }{\partial t}+$

$(-\kappa \left( s\right) u\left( s,t\right) +\tau \left( s\right) w\left(
s,t\right) +\frac{\partial v\left( s,t\right) }{\partial s})\frac{\partial
w\left( s,t\right) }{\partial t})T\left( s\right) +$

$((\tau \left( s\right) v\left( s,t\right) +\frac{\partial w\left(
s,t\right) }{\partial s})\frac{\partial u\left( s,t\right) }{\partial t}-$

$(1-\kappa \left( s\right) v\left( s,t\right) +\frac{\partial u\left(
s,t\right) }{\partial s})\frac{\partial w\left( s,t\right) }{\partial t}%
)N\left( s\right) +$

$((-\kappa \left( s\right) u\left( s,t\right) +\tau \left( s\right) w\left(
s,t\right) +\frac{\partial v\left( s,t\right) }{\partial s})\frac{\partial
u\left( s,t\right) }{\partial t}-$

$(1-\kappa \left( s\right) v\left( s,t\right) +\frac{\partial u\left(
s,t\right) }{\partial s})\frac{\partial v\left( s,t\right) }{\partial t}%
)B\left( s\right) .$

Thus, we get 
\begin{equation*}
n\left( s,t_{0}\right) =\phi _{1}\left( s,t_{0}\right) T\left( s\right)
+\phi _{2}\left( s,t_{0}\right) N\left( s\right) +\phi _{3}\left(
s,t_{0}\right) B\left( s\right) ,
\end{equation*}%
where

$\phi _{1}\left( s,t_{0}\right) =\frac{\partial v\left( s,t_{0}\right) }{%
\partial s}\frac{\partial w\left( s,t_{0}\right) }{\partial t}-\frac{%
\partial w\left( s,t_{0}\right) }{\partial s}\frac{\partial v\left(
s,t_{0}\right) }{\partial t},$

$\phi _{2}\left( s,t_{0}\right) =-\left( 1+\frac{\partial u\left(
s,t_{0}\right) }{\partial s}\right) \frac{\partial w\left( s,t_{0}\right) }{%
\partial t}+\frac{\partial w\left( s,t_{0}\right) }{\partial s}\frac{%
\partial u\left( s,t_{0}\right) }{\partial t},$

$\phi _{3}\left( s,t_{0}\right) =-\left( 1+\frac{\partial u\left(
s,t_{0}\right) }{\partial s}\right) \frac{\partial v\left( s,t_{0}\right) }{%
\partial t}+\frac{\partial v\left( s,t_{0}\right) }{\partial s}\frac{%
\partial u\left( s,t_{0}\right) }{\partial t}.$

This follows that $n_{1}\left( s\right) //n\left( s,t_{0}\right) $, $%
L_{1}\leq s\leq L_{2}$, if and only if there exits a function $\lambda
\left( s\right) \neq 0$ such that%
\begin{equation}
\phi _{1}\left( s,t_{0}\right) =0,~\phi _{2}\left( s,t_{0}\right) =\lambda
\left( s\right) \cosh \theta ,~\phi _{3}\left( s,t_{0}\right) =\lambda
\left( s\right) \sinh \theta .  \tag{4.3}
\end{equation}%
Secondly, since $S\left( T\right) =\omega T,~\omega \neq 0,$%
\begin{equation}
\theta \left( s\right) =-\int\limits_{s_{0}}^{s}\tau ds+\theta _{0\text{,}} 
\tag{4.4}
\end{equation}%
where $s_{0}$ is the starting value of arc length and $\theta =\theta \left(
s\right) .$ In this paper, we assume $s_{0}=0.$

Combining (4.2), (4.3) and (4.4), we have the following theorem.

\begin{theorem}
A spacelike curve $r\left( s\right) $ with timelike binormal is a line of
curvature on the surface $P\left( s,t\right) $ if and only if the followings
are satisfied
\end{theorem}

$\theta \left( s\right) =-\int\limits_{s_{0}}^{s}\tau ds+\theta \left(
0\right) ,$

$u\left( s,t_{0}\right) =v\left( s,t_{0}\right) =w\left( s,t_{0}\right)
\equiv 0,$

$\phi _{1}\left( s,t_{0}\right) \equiv 0~,~\phi _{2}\left( s,t_{0}\right)
=\lambda \left( s\right) \cosh \theta ,~\phi _{3}\left( s,t_{0}\right)
=\lambda \left( s\right) \sinh \theta .$

\bigskip

Now, we analyse two different types of the marching-scale functions $u\left(
s,t\right) ,$

$v\left( s,t\right) ~$and$~w\left( s,t\right) $ in the Eq. (4.1).

\textbf{(i)} If we choose

$u\left( s,t\right) =\overset{p}{\underset{k=1}{\sum }}a_{1k}l\left(
s\right) ^{k}U\left( t\right) ^{k},$ $v\left( s,t\right) =\overset{p}{%
\underset{k=1}{\sum }}a_{2k}m\left( s\right) ^{k}V\left( t\right) ^{k}$ and $%
w\left( s,t\right) =\overset{p}{\underset{k=1}{\sum }}a_{3k}n\left( s\right)
^{k}W\left( t\right) ^{k}$

then, we can simply express the sufficient condition for which the curve $%
r\left( s\right) $ is a line of curvature on the surface $P\left( s,t\right) 
$ as

\begin{eqnarray}
U\left( t_{0}\right) &=&V\left( t_{0}\right) =W\left( t_{0}\right) =0  \notag
\\
\theta \left( s\right) &=&-\int\limits_{s_{0}}^{s}\tau ds+\theta _{0\text{,}}
\TCItag{4.5} \\
a_{21}m\left( s\right) V^{\prime }\left( t_{0}\right) &=&-\lambda \left(
s\right) \sinh \theta ,~a_{31}n\left( s\right) W^{\prime }\left(
t_{0}\right) =-\lambda \left( s\right) \cosh \theta ,  \notag
\end{eqnarray}

\bigskip $\lambda \left( s\right) \neq 0,$ where $l\left( s\right) ,m\left(
s\right) ,n\left( s\right) ,U\left( t\right) ,V\left( t\right) $ and $%
W\left( t\right) $ are $C^{1}$ functions, $a_{ij}\in 
\mathbb{R}
$ $\left( k=1,2,3;j=1,2,3,...,p\right) $.

\textbf{(ii)} If we choose

$u\left( s,t\right) =f\left( \overset{p}{\underset{k=1}{\sum }}a_{1k}l\left(
s\right) ^{k}U\left( t\right) ^{k}\right) ,$ \ $v\left( s,t\right) =g\left( 
\overset{p}{\underset{k=1}{\sum }}a_{2k}m\left( s\right) ^{k}V\left(
t\right) ^{k}\right) $ \ and $\ \ w\left( s,t\right) =h\left( \overset{p}{%
\underset{k=1}{\sum }}a_{3k}n\left( s\right) ^{k}W\left( t\right)
^{k}\right) $

then, we can express the sufficient condition for which the curve $r\left(
s\right) $ is a line of curvature on the surface $P\left( s,t\right) $ as

\begin{eqnarray}
U\left( t_{0}\right) &=&V\left( t_{0}\right) =W\left( t_{0}\right) =0\ and\
f\left( 0\right) =g\left( 0\right) =h\left( 0\right) =0,  \notag \\
\theta \left( s\right) &=&-\int\limits_{s_{0}}^{s}\tau ds+\theta _{0}, 
\TCItag{4.6} \\
g^{^{\prime }}\left( 0\right) a_{21}m\left( s\right) V^{\prime }\left(
t_{0}\right) &=&-\lambda \left( s\right) \sinh \theta ,~h^{^{\prime }}\left(
0\right) a_{31}n\left( s\right) W^{\prime }\left( t_{0}\right) =-\lambda
\left( s\right) \cosh \theta ,  \notag
\end{eqnarray}

$\lambda \left( s\right) \neq 0,$ where $l\left( s\right) ,m\left( s\right)
,n\left( s\right) ,U\left( t\right) ,V\left( t\right) $ and $W\left(
t\right) $ are $C^{1}$ functions, $a_{ij}\in 
\mathbb{R}
$ $\left( k=1,2,3;j=1,2,3,...,p\right) $.

\begin{example}
Let $r\left( s\right) =\left( \cos \left( s\right) ,\sin \left( s\right)
,0\right) $ be an arc-length spacelike curve $0<s\leq 2\pi .$It is easy to
show that
\end{example}

$T\left( s\right) =\left( -\sin \left( s\right) ,\cos \left( s\right)
,0\right) ,$

$N\left( s\right) =\left( -\cos \left( s\right) ,-\sin \left( s\right)
,0\right) ,$

$B\left( s\right) =\left( 0,0,1\right) ,$

$\tau =0.$

If we let $\theta \left( 0\right) =0_{\text{,}}$then we have $\theta \left(
s\right) =0.$ By taking $\lambda \left( s\right) =s,~t_{0}=0$ and the
marching-scale functions as

$u\left( s,t\right) =\sin t,~v\left( s,t\right) =0,~w\left( s,t\right)
=-\sinh \left( ts\right) $

then the Eq. (4.6) is satisfied. So, we have the following surface as a
member of the surface pencil with common line of curvature $r\left( s\right) 
$ (Fig. 4.1) as

$P_{3}\left( s,t\right) =\left( \cos \left( s\right) -\sin \left( t\right)
\sin \left( s\right) ,\sin \left( s\right) +\cos \left( s\right) \sin \left(
t\right) ,\sinh \left( ts\right) \right) ,$

where $0<s\leq 2\pi ,~-1\leq t\leq 1.$

\bigskip

\FRAME{dtbpFU}{2.0911in}{1.8602in}{0pt}{\Qcb{Fig. 4.1. $P_{3}\left(
s,t\right) $ as a member of the surface pencil and its line of curvature}}{}{%
Figure}{\special{language "Scientific Word";type
"GRAPHIC";maintain-aspect-ratio TRUE;display "USEDEF";valid_file "T";width
2.0911in;height 1.8602in;depth 0pt;original-width 2.0522in;original-height
1.823in;cropleft "0";croptop "1";cropright "1";cropbottom "0";tempfilename
'MDS86S02.wmf';tempfile-properties "XPR";}}\bigskip

For the same curve if we choose $t_{0}=0$ , $\lambda \left( s\right) =-\cosh
\left( s\right) $ and the marching-scale functions as

$u\left( s,t\right) =\overset{4}{\underset{k=1}{\sum }}\sin ^{k}\left(
s\right) \sin ^{k}\left( t\right) ,~v\left( s,t\right) =0,~w\left(
s,t\right) =\overset{4}{\underset{k=1}{\sum }}\cosh ^{k}\left( s\right)
\sinh ^{k}\left( t\right) $

we have the Eqn. (4.5) is satisfied. Thus, the surface

$P_{4}\left( s,t\right) =(\cos \left( s\right) -\sin \left( s\right) \overset%
{4}{\underset{k=1}{\sum }}\sin ^{k}\left( s\right) \sin ^{k}\left( t\right)
,\sin \left( s\right) +\cos \left( s\right) \overset{4}{\underset{k=1}{\sum }%
}\sin ^{k}\left( s\right) \sin ^{k}\left( t\right) ,$

$\overset{4}{\underset{k=1}{\sum }}\cosh ^{k}\left( s\right) \sinh
^{k}\left( t\right) )$

is a member of the surface pencil with common line of curvature $r\left(
s\right) $, where $-1.1\leq s\leq 1.1,~-0.7\leq t\leq 0.2$ (Fig.4 2).

\bigskip

\bigskip \FRAME{dtbpFU}{2.7086in}{2.6775in}{0pt}{\Qcb{Fig. 4.2. $P_{4}\left(
s,t\right) $ as a member of the surface pencil and its line of curvature}}{}{%
Figure}{\special{language "Scientific Word";type
"GRAPHIC";maintain-aspect-ratio TRUE;display "USEDEF";valid_file "T";width
2.7086in;height 2.6775in;depth 0pt;original-width 2.6671in;original-height
2.6351in;cropleft "0";croptop "1";cropright "1";cropbottom "0";tempfilename
'MDS86S03.wmf';tempfile-properties "XPR";}}

\bigskip

Now let $P=P\left( s,t\right) $ be a parametric timelike surface and $%
r=r\left( s\right) $ be a timelike curve. The surface is defined by the
given curve as%
\begin{equation}
P\left( s,t\right) =r\left( s\right) +\left( u\left( s,t\right) ,v\left(
s,t\right) ,w\left( s,t\right) \right) \left( 
\begin{array}{c}
T\left( s\right) \\ 
N\left( s\right) \\ 
B\left( s\right)%
\end{array}%
\right) ,  \tag{4.7}
\end{equation}%
$L_{1}\leq s\leq L_{2},T_{1}\leq t\leq T_{2}$.

Let $n_{1}~\left( n_{1}=\cos \theta N+\sin \theta B\right) $ be a vector
orthogonal to the curve $r\left( s\right) $, where $\theta =\theta \left(
s\right) $ is the Lorentzian spacelike angle between $N$ and $n_{1}.$ The
curve $r\left( s\right) $ is a line of curvature on the surface $P\left(
s,t\right) $ if and only if $n_{1}$ is parallel to the normal vector $%
n\left( s,t\right) $ of the surface $P\left( s,t\right) $ and $S\left(
T\right) =\omega T,~\omega \neq 0$, where $S$ is the shape operator of the
surface.

Since the curve $r\left( s\right) $ is a parametric curve on the surface $%
P\left( s,t\right) $, there exists a parameter $t_{0}\in \left[ T_{1},T_{2}%
\right] $ such that $P\left( s,t_{0}\right) =r\left( s\right) $ $L_{1}\leq
s\leq L_{2},$ that is , 
\begin{equation}
u\left( s,t_{0}\right) =v\left( s,t_{0}\right) =w\left( s,t_{0}\right)
\equiv 0~,~L_{1}\leq s\leq L_{2},T_{1}\leq t\leq T_{2}.  \tag{4.8}
\end{equation}%
Fistly, we derive the condition for $n_{1}$ to be parallel to the normal
vector $n\left( s,t\right) $ of the surface $P\left( s,t\right) :$

The normal vector can be expressed as

$n\left( s,t\right) =\frac{\partial P\left( s,t\right) }{\partial s}\times 
\frac{\partial P\left( s,t\right) }{\partial t}$

$=((-\tau \left( s\right) v\left( s,t\right) +\frac{\partial w\left(
s,t\right) }{\partial s})\frac{\partial v\left( s,t\right) }{\partial t}-$

$(\kappa \left( s\right) u\left( s,t\right) +\tau \left( s\right) w\left(
s,t\right) +\frac{\partial v\left( s,t\right) }{\partial s})\frac{\partial
w\left( s,t\right) }{\partial t})T\left( s\right) +$

$((-\tau \left( s\right) v\left( s,t\right) +\frac{\partial w\left(
s,t\right) }{\partial s})\frac{\partial u\left( s,t\right) }{\partial t}-$

$(1+\kappa \left( s\right) v\left( s,t\right) +\frac{\partial u\left(
s,t\right) }{\partial s})\frac{\partial w\left( s,t\right) }{\partial t}%
)N\left( s\right) +$

$(-(\kappa \left( s\right) u\left( s,t\right) +\tau \left( s\right) w\left(
s,t\right) +\frac{\partial v\left( s,t\right) }{\partial s})\frac{\partial
u\left( s,t\right) }{\partial t}+$

$(1+\kappa \left( s\right) v\left( s,t\right) +\frac{\partial u\left(
s,t\right) }{\partial s})(\frac{\partial v\left( s,t\right) }{\partial t}%
))B\left( s\right) $

Thus, we get 
\begin{equation*}
n\left( s,t_{0}\right) =\phi _{1}\left( s,t_{0}\right) T\left( s\right)
+\phi _{2}\left( s,t_{0}\right) N\left( s\right) +\phi _{3}\left(
s,t_{0}\right) B\left( s\right) ,
\end{equation*}%
where

$\phi _{1}\left( s,t_{0}\right) =\frac{\partial w\left( s,t_{0}\right) }{%
\partial s}\frac{\partial v\left( s,t_{0}\right) }{\partial t}-\frac{%
\partial v\left( s,t_{0}\right) }{\partial s}\frac{\partial w\left(
s,t_{0}\right) }{\partial t},$

$\phi _{2}\left( s,t_{0}\right) =\frac{\partial w\left( s,t_{0}\right) }{%
\partial s}\frac{\partial u\left( s,t_{0}\right) }{\partial t}-\left( 1+%
\frac{\partial u\left( s,t_{0}\right) }{\partial s}\right) \frac{\partial
w\left( s,t_{0}\right) }{\partial t},$

$\phi _{3}\left( s,t_{0}\right) =\left( 1+\frac{\partial u\left(
s,t_{0}\right) }{\partial s}\right) \frac{\partial v\left( s,t_{0}\right) }{%
\partial t}-\frac{\partial v\left( s,t_{0}\right) }{\partial s}\frac{%
\partial u\left( s,t_{0}\right) }{\partial t}.$

This follows that $n_{1}\left( s\right) //n\left( s,t_{0}\right) $, $%
L_{1}\leq s\leq L_{2}$, if and only if there exits a function $\lambda
\left( s\right) \neq 0$ such that%
\begin{equation}
\phi _{1}\left( s,t_{0}\right) =0,~\phi _{2}\left( s,t_{0}\right) =\lambda
\left( s\right) \cos \theta ,~\phi _{3}\left( s,t_{0}\right) =\lambda \left(
s\right) \sin \theta .  \tag{4.9}
\end{equation}%
Secondly, since $S\left( T\right) =\omega T,~\omega \neq 0,$%
\begin{equation}
\theta \left( s\right) =\int\limits_{s_{0}}^{s}\tau ds+\theta _{0\text{,}} 
\tag{4.10}
\end{equation}%
where $s_{0}$ is the starting value of arc length and $\theta =\theta \left(
s\right) .$ In this paper, we assume that $s_{0}=0.$

Combining (4.8), (4.9) and (4.10), we have the following theorem.

\begin{theorem}
A timelike curve $r\left( s\right) $ is a line of curvature on the surface $%
P\left( s,t\right) $ if and only if the followings are satisfied
\end{theorem}

$\theta \left( s\right) =\int\limits_{s_{0}}^{s}\tau ds+\theta \left(
0\right) ,$

$u\left( s,t_{0}\right) =v\left( s,t_{0}\right) =w\left( s,t_{0}\right)
\equiv 0,$

$\phi _{1}\left( s,t_{0}\right) \equiv 0~,~\phi _{2}\left( s,t_{0}\right)
=\lambda \left( s\right) \cos \theta ,~\phi _{3}\left( s,t_{0}\right)
=\lambda \left( s\right) \sin \theta .$

\bigskip

We call the set of surfaces defined by (4.1) - (4.4) or (4-7) - (4-10) 
\textit{timelike surface pencil with a common line of curvature}. Any
surface $P\left( s,t\right) $ satisfying these conditions is a member of
this family.

\bigskip

Now, we analyse two different types of the marching-scale functions

$u\left( s,t\right) ,~v\left( s,t\right) ~$and$~w\left( s,t\right) $ in the
Eq. (4.1).

\textbf{(i)} If we choose marching-scale functions as

$u\left( s,t\right) =\overset{p}{\underset{k=1}{\sum }}a_{1k}l\left(
s\right) ^{k}U\left( t\right) ^{k},$ $v\left( s,t\right) =\overset{p}{%
\underset{k=1}{\sum }}a_{2k}m\left( s\right) ^{k}V\left( t\right) ^{k}$ and $%
w\left( s,t\right) =\overset{p}{\underset{k=1}{\sum }}a_{3k}n\left( s\right)
^{k}W\left( t\right) ^{k}$

then, we can simply express the sufficient condition for which the curve $%
r\left( s\right) $ is a line of curvature of the surface $P\left( s,t\right) 
$ as

\begin{eqnarray}
U\left( t_{0}\right) &=&V\left( t_{0}\right) =W\left( t_{0}\right) =0, 
\notag \\
\theta \left( s\right) &=&\int\limits_{s_{0}}^{s}\tau ds+\theta , 
\TCItag{4.11} \\
a_{21}m\left( s\right) V^{\prime }\left( t_{0}\right) &=&\lambda \left(
s\right) \sin \theta ,~a_{31}n\left( s\right) W^{\prime }\left( t_{0}\right)
=-\lambda \left( s\right) \cos \theta ,  \notag
\end{eqnarray}

\bigskip $\lambda \left( s\right) \neq 0,$ where $l\left( s\right) ,m\left(
s\right) ,n\left( s\right) ,U\left( t\right) ,V\left( t\right) $ and $%
W\left( t\right) $ are $C^{1}$ functions, $a_{ij}\in 
\mathbb{R}
$ $\left( k=1,2,3;j=1,2,3,...,p\right) .$

\textbf{(ii)} If we choose marching-scale functions as

$u\left( s,t\right) =f\left( \overset{p}{\underset{k=1}{\sum }}a_{1k}l\left(
s\right) ^{k}U\left( t\right) ^{k}\right) ,$ \ $v\left( s,t\right) =g\left( 
\overset{p}{\underset{k=1}{\sum }}a_{2k}m\left( s\right) ^{k}V\left(
t\right) ^{k}\right) $ \ and $\ \ w\left( s,t\right) =h\left( \overset{p}{%
\underset{k=1}{\sum }}a_{3k}n\left( s\right) ^{k}W\left( t\right)
^{k}\right) $

then, we can express the sufficient condition for which the curve $r\left(
s\right) $ is a line of curvature on the surface $P\left( s,t\right) $ as

\begin{eqnarray}
U\left( t_{0}\right) &=&V\left( t_{0}\right) =W\left( t_{0}\right) =0\ and\
f\left( 0\right) =g\left( 0\right) =h\left( 0\right) =0,  \notag \\
\theta \left( s\right) &=&\int\limits_{s_{0}}^{s}\tau ds+\theta _{0\text{,}}
\TCItag{4.12} \\
g^{^{\prime }}\left( 0\right) a_{21}m\left( s\right) V^{\prime }\left(
t_{0}\right) &=&\lambda \left( s\right) \sin \theta ,~h^{^{\prime }}\left(
0\right) a_{31}n\left( s\right) W^{\prime }\left( t_{0}\right) =-\lambda
\left( s\right) \cos \theta ,  \notag
\end{eqnarray}

$\lambda \left( s\right) \neq 0,~$\ where $l\left( s\right) ,m\left(
s\right) ,n\left( s\right) ,U\left( t\right) ,V\left( t\right) $ and $%
W\left( t\right) $ are $C^{1}$ functions, $a_{ij}\in 
\mathbb{R}
$ $\left( k=1,2,3;j=1,2,3,...,p\right) $.

\begin{example}
Let $r\left( s\right) =\left( \cosh \left( s\right) ,0,\sinh \left( s\right)
\right) $be an arc-length timelike curve ,$0$ $\leq s\leq 2\pi .$It is easy
to show that
\end{example}

$T\left( s\right) =\left( \sinh \left( s\right) ,0,\cosh \left( s\right)
\right) ,$

$N\left( s\right) =\left( \cosh \left( s\right) ,0,\sinh \left( s\right)
\right) ,$

$B\left( s\right) =\left( 0,-1,0\right) ,$

$\tau =0.$

By taking $\lambda \left( s\right) =-s,~t_{0}=0$ and the marching-scale
functions as

$u\left( s,t\right) =\sinh \left( t\right) ,~v\left( s,t\right) =0,~w\left(
s,t\right) =\cosh \left( st\right) $

we have the Eq. (4.12) is satisfied. So, we obtain the following surface as
a member of the surface pencil with common line of curvature $r\left(
s\right) $ (Fig. 4.3) as

$P_{5}\left( s,t\right) =\left( \cosh \left( s\right) +\sinh \left( t\right)
\sinh \left( s\right) ,-\cosh \left( s\ast t\right) ,\sinh \left( s\right)
+\cosh \left( s\right) \sinh \left( t\right) \right) ,$

where $0$ $\leq s\leq 2\pi ,$ $-1\leq t\leq 1.$

\bigskip \FRAME{dtbpFU}{3.557in}{3.4627in}{0pt}{\Qcb{Fig. 4.3. $P_{5}\left(
s,t\right) $ as a member of the surface pencil and its line of curvature}}{}{%
Figure}{\special{language "Scientific Word";type
"GRAPHIC";maintain-aspect-ratio TRUE;display "USEDEF";valid_file "T";width
3.557in;height 3.4627in;depth 0pt;original-width 3.5103in;original-height
3.4169in;cropleft "0";croptop "1";cropright "1";cropbottom "0";tempfilename
'MDS86S04.wmf';tempfile-properties "XPR";}}

For the same curve let $\lambda \left( s\right) =-\sinh \left( s\right)
,~t_{0}=0$ and the marching-scale functions

$u\left( s,t\right) =\overset{4}{\underset{k=1}{\sum }}\sinh ^{k}\left(
t\right) ,~v\left( s,t\right) \equiv 0,~w\left( s,t\right) =\overset{4}{%
\underset{k=1}{\sum }}\sinh ^{k}\left( s\right) \sinh ^{k}\left( t\right) .$

Now we have the Eq. (4.11) is satisfied. Thus, the surface

$P_{6}\left( s,t\right) =(\cosh \left( s\right) +\sinh \left( s\right) 
\overset{4}{\underset{k=1}{\sum }}\sinh ^{k}\left( t\right) ,-\overset{4}{%
\underset{k=1}{\sum }}\sinh ^{k}\left( s\right) \sinh ^{k}\left( t\right)
,\sinh \left( s\right) +\cosh \left( s\right) \overset{4}{\underset{k=1}{%
\sum }}\sinh ^{k}\left( t\right) )$

is a member of the surface pencil with common line of curvature $r\left(
s\right) $, where $-1$ $\leq s\leq 1,$ $-0.4\leq t\leq 0.4.$

\FRAME{dtbpFU}{2.4258in}{2.6351in}{0pt}{\Qcb{Fig. 4.4. $P_{6}\left(
s,t\right) $ as a member of the surface pencil and its line of curvature}}{}{%
Figure}{\special{language "Scientific Word";type
"GRAPHIC";maintain-aspect-ratio TRUE;display "USEDEF";valid_file "T";width
2.4258in;height 2.6351in;depth 0pt;original-width 2.3851in;original-height
2.5936in;cropleft "0";croptop "1";cropright "1";cropbottom "0";tempfilename
'MDS86S05.wmf';tempfile-properties "XPR";}}

\begin{acknowledgement}
The authors appreciate the comments and valuable suggestions of the editor
and the reviewer. Their advice helped to improve the clarity and
presentation of this paper. The second author would like to thank TUBITAK
(The Scientific and Technological Research Council of Turkey) for their
financial supports during his doctorate studies.
\end{acknowledgement}

\section{\protect\bigskip References}

\ \ \ \ \ [1] Wang, G.J., Tang, K., Tai, C.L., Parametric representantion of
a surface pencil with a common spatial geodesic, Computer Aided Design,
2004, 36:\textbf{\ }447-459.

[2] Kasap, E., Aky\i ld\i z, F.T., Surfaces with common geodesic in
Minkowski 3-space, Applied Mathematics and Computation, 2008, 177: 260-270.

[3] Li, C.Y., Wang, R. H., Zhu, C.G., Parametric representation of a surface
pencil with a common line of curvature, Computer Aided Design, 2011, 43:
1110-1117.

[4] \c{S}affak, G., Kasap, E., Family of surface with a common null
geodesic, International Journal of Physical Sciences, 2009, Vol.4 (8):
428-433.

[5] O'Neill, B., Semi-Riemannian Geometry with applications to relativity,
New York: Academic Press, 1983.

[6] Akutagawa, K. Nishikawa, S., The Gauss map and spacelike surface with
prescribed mean curvature in Minkowski 3-space, Tohoku Math. J., 1990, 42:
67-82.

[7] Woestijne, V.D.I, Minimal Surface of the 3-Dimensional Minkowski Space,
Singapore: World Scientific Publishing, 1990.

[8] Ratcliffe, J.G., Foundations of Hyperbolic Manifolds, New York:
Springer-Verlag, 1994.

[9] Beem, J.K., Ehrlich, P.E., Global Lorentzian Geometry, New York: Marcel
Dekker, 1981.$\ \ $

\end{document}